\documentclass[12pt]{article}
\usepackage{latexsym,amsfonts,amssymb,amsmath}
\usepackage[dvips]{color}
\usepackage{colordvi,multicol}

\numberwithin{equation}{section}

\def \cal{\mathcal}

\newtheorem{thm}{Theorem}[section]

\newtheorem{lem}[thm]{Lemma}

\newtheorem{rem}[thm]{Remark}

\date{}
\begin{document}

\title{\bf Some explorations on two conjectures about Rademacher sequences}
 \author{ Ze-Chun Hu, Guolie Lan\footnote{Corresponding author}, Wei Sun\\ \\
 {\small College of Mathematics, Sichuan  University,  China}\\
 {\small zchu@scu.edu.cn}\\ \\
 {\small School of Economics and Statistics, Guangzhou University, China}\\
 {\small langl@gzhu.edu.cn}\\ \\
 {\small Department of Mathematics and Statistics, Concordia
University, Canada}\\
{\small wei.sun@concordia.ca}}

\maketitle

\begin{abstract}
\noindent In this paper, we explore two conjectures about Rademacher sequences. Let $(\epsilon_i)$ be a Rademacher sequence, i.e., a sequence of independent  $\{-1,1\}$-valued symmetric random variables. Set  $S_n=a_1\epsilon_1+\cdots+a_n\epsilon_n$ for $a=(a_1,\dots,a_n)\in \mathbb{R}^n$.
The first conjecture says that $P\left(\left|S_n\right|\leq \|a\|\right)\geq\frac{1}{2}$ for all $a\in \mathbb{R}^n$ and $n\in \mathbb{N}$. The second conjecture says that $P\left(\left|S_n\right|\geq\|a\|\right)\geq \frac{7}{32}$ for all $a\in \mathbb{R}^n$ and $n\in \mathbb{N}$.
Regarding the first conjecture, we present several new equivalent formulations. These include a topological view, a combinatorial version and a strengthened version of the conjecture. Regarding the second conjecture, we prove that
it holds true when $n\leq 7$.

\noindent

\end{abstract}




\section{Introduction}\setcounter{equation}{0}

Let $a=(a_1,\dots,a_n)$ be a real vector with $\sum_{i=1}^na_i^2=1$.
Consider the $2^n$ sign sums of the form $S=\pm a_1\pm\cdots\pm a_n$.
In 1986, B. Tomaszewaki (see Guy \cite{Gu86}) posed the question that
whether at least half of these sign sums satisfy $|S|\leq 1$.
This question looks simple but turns out to be very challenging.
Although several partial results have been obtained,
the conjecture still remains open after more than 30 years.

Tomaszewaki's conjecture can  also be described as a probability problem.
Throughout this paper, we let $\{\epsilon_i,i=1,2,\ldots\}$ be a Rademacher sequence, i.e., a sequence of independent random variables with  $P(\epsilon_i=-1)=P(\epsilon_i=1)=\frac{1}{2}$.
Let $n\in \mathbb{N}$ and $\epsilon=(\epsilon_1,\dots,\epsilon_n)$.
For $a=(a_1,\dots,a_n)\in \mathbb{R}^n$, denote $a\cdot\epsilon:=a_1\epsilon_1+\cdots+a_n\epsilon_n$.  Then Tomaszewaki's conjecture can be formulated as
$$
 P\left(\left|a\cdot\epsilon\right|\leq\|a\|\right)\geq\frac{1}{2}.
$$
Hereafter $\|a\|:=(\sum_{i=1}^n{a_i}^2)^{\frac{1}{2}}$.  Consider the non-increasing sequence of constants
    $$T_n=\inf\left\{ P\left(\left|a\cdot\epsilon\right|\leq\|a\|\right):{a\in \mathbb{R}^n}\right\}, \quad n\in\mathbb{N}.$$
$T=\inf_n \{T_n\}$ is known as  Tomaszewaki's constant. This constant occurs frequently in optimization and operations research
(cf. \cite{DDS16,To18}).
Since $T_2=\frac{1}{2}$  whenever $a_1a_2\neq0$,
Tomaszewaki's conjecture is equivalent to $T=\frac{1}{2}$.

Holzman and Kleitman \cite{HK92} proved that $T\geq \frac{3}{8}$.
Ben-Tal, Nemirovski and Roos \cite{BNR02} independently proved that $T\geq \frac{1}{3}$
(see also \cite{Sh12} for an improvement to $36\%$).
The update results are $T>\frac{13}{32}$ by Boppana and Holzman \cite{BH17}
and $T_9=\frac12$ by Hendriks and van Zuijlen \cite{HZ17}.
van Zurjlen \cite{Zu11} and von Heymann \cite{He12} showed that Tomaszewaki's conjecture is true for the special case that all $|a_i|$ are equal.
Tomaszewaki's constant can be approximated surprisingly via an algorithm constructed by De et al. \cite{DDS16}. According to that algorithm, for any $\varepsilon>0$, a value between $T$ and $T+\varepsilon$ is returned when the running time is long enough.
In this paper, we will present several new equivalent formulations for Tomaszewaki's conjecture. These include a topological view, a combinatorial version and a strengthened version of the conjecture. We hope that our results can shed new light on this long-standing open problem.

Besides $T$, another constant related to the Rademacher sequence arises in a variety of contexts.
Consider the non-increasing sequence of constants
\begin{equation}\label{Gn}
G_n=\inf\left\{ P\left(|a\cdot\epsilon|\geq\|a\|\right):{a\in \mathbb{R}^n}\right\}, \quad n\in\mathbb{N}.
\end{equation}
$G=\inf_n \{G_n\}$ is known as Hitczenko and Kwapie\'{n}'s
constant. This  constant has been applied to Kolmogorov-type lower estimates for independent random variables
and other problems (cf. \cite{Bu67,HK94}).
Burkholder \cite{Bu67} showed that $G>0$. Hitczenko and Kwapie\'{n} \cite{HK94} showed that $G>\frac{1}{4} e^{-4}$. Oleszkiewicz \cite{Ol96} showed that $G>\frac{1}{10}$. He also pointed out that $G$ might be equal to  $\frac{7}{32}$ and proved that $G_6=\frac{7}{32}$.
It can be checked that $G_1=1,\ G_2=\frac{1}{2},\ G_3=G_4=G_5=\frac{1}{4}$.
In this paper, we will prove that  $G_7=\frac{7}{32}$ (see Theorem \ref{thm-21} and Remark \ref{rem-21} below).

Note that $G_n$ yields an upper bound $1-G_n$ for
$P\left(\left|a\cdot\epsilon\right|<\|a\|\right)$,
but not for $P\left(\left|a\cdot\epsilon\right|\leq\|a\|\right)$.
The following slightly different constants yield an upper bound for $P\left(\left|a\cdot\epsilon\right|\leq\|a\|\right)$:
\begin{equation}\label{Gn'}
G'_n:=\inf\left\{ P\left(|a\cdot\epsilon|>\|a\|\right):{a\in \mathbb{R}^n}, a_i\neq0,\forall i=1,\ldots,n\right\}.
\end{equation}
We will show that $
G_1'=0,\ G_2'=\frac{1}{2},\ G_3'=G_5'=\frac{1}{4},\ G_4'=\frac{1}{8},\ G_6'=\frac{3}{16},\ G_7'=\frac{7}{32}.
$

The rest of this paper is organized as follows.
In Section 2, we explore Tomaszewaki's conjecture and
present several new equivalent formulations.
In Section 3, we study Hitczenko and Kwapie\'{n}'s conjecture and prove among other things that $G_7=G_7'=\frac{7}{32}$.

\section{Equivalent formulations for Tomaszewski's conjecture}\setcounter{equation}{0}

In this section, we make some explorations on Tomaszewski's conjecture. In Section 2.1, we give an equivalent formulation for Tomaszewski's conjecture from a topological view. In Section 2.2, we reformulate Tomaszewski's conjecture as a conjecture in terms of natural numbers. In Section 2.3, we present a strengthened version of Tomaszewski's conjecture.

\subsection{A topological view}

Suppose $n\ge 2$. Define
$$
{\cal S}^{n-1}:=\left\{a\in \mathbb{R}^n:\|a\|=1\right\},
$$
and
$$
{\cal D}:=\left\{a\in {\cal S}^{n-1}: P\left(\left|a\cdot\epsilon\right|\le1\right)\ge \frac{1}{2}\right\}.
$$

\begin{lem}\label{pro-2.1}
Tomaszewski's conjecture holds if and only if there exists a dense subset $V$ of ${\cal S}^{n-1}$ such that the inequality
\begin{equation}\label{aforTC}
P\left(\left|a\cdot\epsilon\right|\le1\right)\ge \frac{1}{2}
\end{equation}
 holds for all $a\in V$.
\end{lem}
{\bf Proof.} We need only prove the sufficiency. To this end, it is sufficient to show that ${\cal D}$ is a closed subset of ${\cal S}^{n-1}$.

Let $a\in {\cal S}^{n-1}$ and $\{a^{(m)}\}\subset {\cal D}$ such that $\|a^{(m)}-a\|\rightarrow 0$ as $m\rightarrow\infty$.  We have
$$
P(|a^{(m)}\cdot\epsilon|\le 1)\ge \frac{1}{2},\ \ \forall m\in \mathbb{N}.
$$
Note that $\{-1,1\}^n$ contains $2^n$ vectors. Then, there exist
$$\{(\eta^{(l)}_1,\dots,\eta^{(l)}_n)\in\{-1,1\}^n:1\le l\le 2^{n-1}\}
$$
and a subsequence $\{m_k\}$ such that
$$
\left|\sum_{i=1}^na^{(m_k)}_i\eta^{(l)}_i\right|\le1,\ \ \forall l\in \{1,2,\ldots,2^{n-1}\},\ \forall k\in\mathbb{N}.
$$

By $\|a^{(m)}-a\|\rightarrow 0$ as $m\rightarrow\infty$, we get
$$
\left|\sum_{i=1}^na_i\eta^{(l)}_i\right|\le1,\ \ \forall l\in \{ 1,2,\ldots, 2^{n-1}\}.
$$
It follows that
$
    P(|a\cdot\epsilon|\le1)\ge \frac{1}{2},
$
which implies that $a\in {\cal D}$. Hence ${\cal D}$ is a closed subset of ${\cal S}^{n-1}$. The proof is complete. \hfill\fbox

To make further analysis, we define
\begin{eqnarray*}
&&A:=\{a\in {\cal S}^{n-1}:P(|a\cdot\epsilon|=1)=0\},\\
&&B:=\{a\in {\cal S}^{n-1}:P(|a\cdot\epsilon|=1)>0\},
\end{eqnarray*}
and
$$
C:=\left\{a\in A: P(|a\cdot\epsilon|>1)\le \frac{1}{2}\right\}.
$$
Note that
$$
C=\left\{a\in A: P(|a\cdot\epsilon|\ge1)\le \frac{1}{2}\right\}=\left\{a\in A: P(|a\cdot\epsilon|\le1)\ge \frac{1}{2}\right\}.
$$

\begin{lem}\label{lem-2.3}
$A$ is a dense open subset of ${\cal S}^{n-1}$ and $C$ is both an open and closed subset of $A$.
\end{lem}
{\bf Proof.} Obviously, $A$ is a dense subset of ${\cal S}^{n-1}$. We first show that $B$ is a closed subset of $\mathbb{R}^n$ and hence $A={\cal S}^{n-1}- B$ is an open subset of ${\cal S}^{n-1}$. In fact, let $\{a^{(m)}\}\subset B$ and $\|a^{(m)}-a\|\rightarrow 0$ as $m\rightarrow\infty$. Then, there exist $(\eta_1,\dots,\eta_n)\in\{-1,1\}^n$ and a subsequence $\{m_k\}$ such that
$$
\left|\sum_{i=1}^na^{(m_k)}_i\eta_i\right|=1,\ \ \forall k\in \mathbb{N}.
$$
By $\|a^{(m)}-a\|\rightarrow 0$ as $m\rightarrow\infty$, we get
$
\left|\sum_{i=1}^na_i\eta_i\right|=1,
$
which implies that $a\in B$. Hence $B$ is a closed subset of $\mathbb{R}^n$.

Next we show that $C$ is an open subset of $A$. Let $a\in C$. Then, there exist
$$\{(\eta^p_1,\dots,\eta^p_n)\in\{-1,1\}^n:1\le p\le 2^{n-1}\}
$$ such that
$$
\left|\sum_{i=1}^na_i\eta^p_i\right|<1,\ \ \forall p\in\{1,2,\ldots,2^{n-1}\}.
$$
Hence there exists $\delta>0$ such that for any $b\in A$ satisfying $\|b-a\|<\delta$,
$$
\left|\sum_{i=1}^nb_i\eta^p_i\right|<1,\ \ \forall p\in \{1,2,\ldots, 2^{n-1}\}.
$$
Thus
$$
P(|b\cdot\epsilon|\ge1)\le \frac{1}{2},
$$
which implies that $b\in C$. Since $a\in C$ is arbitrary, $C$ is an open subset of $A$.

Finally, we show that $C$ is a closed subset of $A$. Let $a\in A$ and $\{a^{(m)}\}\subset C$ such that $\|a^{(m)}-a\|\rightarrow 0$ as $m\rightarrow\infty$. We have
$$
P(|a^{(m)}\cdot\epsilon|\le 1)\ge \frac{1}{2},\ \ \forall m\in \mathbb{N}.
$$
Then, there exist
$$\{(\eta^{(q)}_1,\dots,\eta^{(q)}_n)\in\{1,-1\}^n:1\le q\le 2^{n-1}\}
$$
and a subsequence $\{m_k\}$ such that
$$
\left|\sum_{i=1}^na^{(m_k)}_i\eta^{(q)}_i\right|\le1,\ \ \forall q\in \{1,2,\ldots,2^{n-1}\},\ \forall k\in\mathbb{N}.
$$
By $\|a^{(m)}-a\|\rightarrow 0$ as $m\rightarrow\infty$, we get
$$
\left|\sum_{i=1}^na_i\eta^{(q)}_i\right|\le1,\ \ \forall q\in \{1,2,\ldots,2^{n-1}\}.
$$
Thus $P(|a\cdot\epsilon|\le1)\ge \frac{1}{2},$ which implies that  $a\in C$. Hence $C$ is a closed subset of $A$.
\hfill\fbox

\vskip 0.3cm
By Lemma \ref{lem-2.3}, $A$ is a open subset of ${\cal S}^{n-1}$. We equip $A$ with the inherited topology from ${\cal S}^{n-1}$. Now we can state the main result of this subsection.
\begin{thm}\label{pro-2.2}
Tomaszewski's conjecture holds if and only if each connected component of $A$ contains
at least one point that satisfies (\ref{aforTC}).
\end{thm}
{\bf Proof.} We need only  prove the sufficiency. Suppose that each connected component of $A$ contains
at least one point that satisfies (\ref{aforTC}). Let $F$ be a connected component of $A$. Then there exists a point $a\in F$ satisfying (\ref{aforTC}). By Lemma \ref{lem-2.3}, all points in $F$ must satisfy (\ref{aforTC}). Since $F$ is arbitrary, we conclude that all points in $A$ satisfy (\ref{aforTC}). Finally, since $A$ is a dense subset of ${\cal S}^{n-1}$, Lemma \ref{pro-2.1} implies that Tomaszewski's conjecture holds. \hfill\fbox

\subsection{A combinatorial version}

Let $n\in\mathbb{N}$ and $I_n=\{1,\dots,n\}$. For ${J}\subseteq I_n$, set
${J}_c=I_n-{J}$ and
\begin{equation}\label{JJ}
J^{(2)}=\{(i,j): \ i,j\in J, \ i<j\}.
\end{equation}
For a finite set $A$, we use $|A|$ to denote its cardinality.

\begin{thm}\label{thm-2.4}
Tomaszewski's conjecture holds if and only if for any $ n\in \mathbb{N}$ and $ {l}_1,\dots,{l}_n\in\mathbb{N}$,
\begin{equation}\label{eqn334}
    \frac{1}{2^n}\left|\left\{J\subseteq I_n:\sum_{(i,j)\in J^{(2)}\cup J_c^{(2)}} {l}_i{l}_j-
    \sum_{(i,j)\in J\times J_c} {l}_i{l}_j\leq 0\right\}\right|\geq\frac{1}{2}.
\end{equation}
\end{thm}

Before proving Theorem \ref{thm-2.4}, we present a lemma.

\begin{lem}\label{lem-2.5}
Let $n\geq 2$. Then the following two statements are equivalent:

\noindent
(i) For any $a=(a_1,\dots,a_n)\in\mathbb{R}^n$ and $\epsilon=(\epsilon_1,\dots,\epsilon_n)$,
\begin{equation}\label{eqn3-1}
    P\left(\left|a\cdot\epsilon\right|\leq\|a\|\right)\geq\frac{1}{2}.
\end{equation}
(ii) For any independent random variables $\eta_1,\dots,\eta_n$ with symmetric two-point
distributions,
\begin{equation}\label{eqn3-4}
P\left(\sum_{i\not= j}\eta_i\eta_j\le 0\right)\ge\frac{1}{2}.
\end{equation}
\end{lem}

\noindent {\bf Proof.} (i) $\Rightarrow$ (ii): Let $\eta_1,\dots,\eta_n$ be independent random
variables with symmetric two-point distributions. Define
$$
a_i={|\eta_i|},\quad \epsilon_i=\frac{\eta_i}{|\eta_i|},\quad 1\le i\le n.
$$
Then $a_1,\dots,a_n$ are constants, $\epsilon_1,\dots,\epsilon_n$
are independent  Rademacher random variables with
$$
\eta_i=a_i\epsilon_i,\quad 1\le i\le n.
$$
Note that
\begin{eqnarray}\label{wrong11}
|a\cdot\epsilon|\leq\|a\|
&\Leftrightarrow&\left(\sum_{i=1}^na_i\epsilon_i\right)^2\le\sum_{i=1}^na_i^2\nonumber\\
&\Leftrightarrow&\sum_{i\not=
j}a_i\epsilon_ia_j\epsilon_j\le0\nonumber\\
&\Leftrightarrow&\sum_{i\not= j}\eta_i\eta_j\le0.
\end{eqnarray}
Hence (i) implies (ii).

(ii) $\Rightarrow$ (i): Let $a=(a_1,\dots,a_n)\in\mathbb{R}^n$ and $\epsilon=(\epsilon_1,\dots,\epsilon_n)$. Define $\eta_i=a_i\epsilon_i,1\le i\le n$. Then the proof is complete by (\ref{wrong11}).\hfill\fbox

\vskip 0.3cm
\noindent {\bf Proof of Theorem \ref{thm-2.4}.} ``$\Rightarrow$" Suppose that Tomaszewski's conjecture holds. Let ${l}_1,\dots,{l}_n\in\mathbb{N}$ and
$\epsilon_1,\dots,\epsilon_n$ be independent Rademacher random variables. Define
$$
\eta_i=l_i\varepsilon_i,\ \ 1\le i\le n.
$$
Then, $\eta_1,\dots,\eta_n$ are independent and have symmetric
two-point distributions. By (\ref{eqn3-4}), we get
$$
P\left(\sum_{i\not= j}l_il_j\cdot\varepsilon_i\varepsilon_j\le 0\right)\ge\frac{1}{2}.
$$
Hence there exist $2^{n-1}$ different vectors 
$$s^{(k)}=(s_{k,1},\dots,s_{k,n})\in\{-1,1\}^n, \quad k=1,\dots, 2^{n-1},$$
such that
\begin{equation}\label{eqn337}
    \sum_{i\not= j}l_il_j\cdot s_{k,i}s_{k,j}\le 0, \quad \forall k\in \{1,2,\ldots, 2^{n-1}\}.
\end{equation}

For each vector $s^{(k)}$, define
\begin{equation}\label{eqn33Jk}
    J_k=\{i:s_{k,i}=1\}, \quad k=1,\dots, 2^{n-1}.
\end{equation}
Then $J_k, \ k=1,\dots, 2^{n-1}$ are different subsets of $I_n$.
Note that $s_{k,i}s_{k,j}=1$ for $(i,j)\in  J_k^{(2)}\cup
(J_k)_c^{(2)}$, and $s_{k,i}s_{k,j}=-1$ for $(i,j)\in J_k\times
(J_k)_c$ . It follows from (\ref{eqn337}) that
\begin{equation*}
    \sum_{(i,j)\in J_k^{(2)}\cup (J_k)_c^{(2)}} {l}_i{l}_j-\sum_{(i,j)\in J_k\times (J_k)_c} {l}_i{l}_j\leq 0,
    \quad \forall k\in \{1,\dots, 2^{n-1}\}.
\end{equation*}
Thus (\ref{eqn334}) holds.

 ``$\Leftarrow$"  Suppose that for any $ n\in \mathbb{N}$ and $ {l}_1,\dots,{l}_n\in\mathbb{N}$, \eqref{eqn334} holds. We will show that for any $n\in \mathbb{N}$ and $a=(a_1,\dots,a_n)\in\mathbb{R}^n$, \eqref{eqn3-1} holds.

 We use proof by contradiction.  Suppose that
(\ref{eqn3-1}) does not hold for some $n\in \mathbb{N}$ and
$a=(a_1, \dots, a_n)\in \mathbb{R}^n$. We can assume without loss
of generality that $n\ge 2$ and each $a_i>0$. Then there exist
$2^{n-1}+1$ different vectors
$s^{(k)}=(s_{k,1},\dots,s_{k,n})\in\{-1,1\}^n$, $k=1,\dots,
2^{n-1}+1$ such that $|a\cdot s^{(k)}|>\|a\|$, i.e.,
\begin{equation}\label{eqn3341}
    (a\cdot s^{(k)})^2-a\cdot a>0, \quad \forall k\in \{1,\dots, 2^{n-1}+1\}.
\end{equation}

For each $1\le i\le n$, we choose a sequence
$\{r^{(m)}_i:m=1,2,\dots\}$ of positive rational numbers such that
$r^{(m)}_i\rightarrow a_i$ as $m\rightarrow\infty$. Let
 $$
 r^{(m)}=(r^{(m)}_1,\dots,r^{(m)}_n).
 $$
Then we have $\|r^{(m)}-a\|\rightarrow0$ as $m\rightarrow\infty$.
By (\ref{eqn3341}) and the continuity of functions
$\varphi_k(\cdot): \varphi_k(x)=(x\cdot s^{(k)})^2-x\cdot x$,
$x\in\mathbb{R}^n, k=1,\dots, 2^{n-1}+1$, there exists $m_0\in
\mathbb{N}$ such that
\begin{equation}\label{eqn3342}
    (r^{(m_0)}\cdot s^{(k)})^2-r^{(m_0)}\cdot r^{(m_0)}>0, \quad \forall k\in \{1,\dots, 2^{n-1}+1\}.
\end{equation}

Note that (\ref{eqn3342}) also holds with $r^{(m_0)}$ replaced by
$c r^{(m_0)}$ for any $c>0$. Since $r^{(m_0)}$ is a vector made up
of positive rational numbers, there exists a positive integer $l_0$
such that $l_i=l_0 r^{(m_0)}_i$ are positive integers for
$i=1,\dots,n$. Set $ l=(l_1,\dots,l_n)=l_0 r^{(m_0)}$. Then, $l\in
\mathbb{N}^n$ and
$$
(l\cdot s^{(k)})^2-l\cdot l>0, \quad \forall k\in \{1,\dots, 2^{n-1}+1\}.
$$
Since $s^{(k)}=(s_{k,1},\dots,s_{k,n})\in\{-1,1\}^n$, the above
equalities can be written as
\begin{equation}\label{erg}
    \sum_{i\not= j}l_il_j\cdot s_{k,i}s_{k,j}>0, \quad \forall k\in \{1,\dots, 2^{n-1}+1\}.
\end{equation}

Let $J_k$ be defined by (\ref{eqn33Jk}). Then $s_{k,i}s_{k,j}=1$
for $(i,j)\in  J_k^{(2)}\cup (J_k)_c^{(2)}$ and
$s_{k,i}s_{k,j}=-1$ for $(i,j)\in J_k\times (J_k)_c$ . Thus
(\ref{erg}) implies that
\begin{equation*}
  \sum_{(i,j)\in J_k^{(2)}\cup (J_k)_c^{(2)}} {l}_i{l}_j-\sum_{(i,j)\in J_k\times (J_k)_c} {l}_i{l}_j>0
\end{equation*}
for $2^{n-1}+1$ different subsets $J_k$, $k=1,\dots, 2^{n-1}+1$.
Hence (\ref{eqn334}) does not hold. We have arrived at a contradiction. Therefore, for any $a=(a_1,\dots,a_n)\in\mathbb{R}^n$, \eqref{eqn3-1} holds, i.e., Tomaszewski's conjecture holds. The proof is complete.  \hfill\fbox

\subsection{A strengthened version}

Tomaszewski's conjecture says that for any $n\in \mathbb{N}$ and $a=(a_1,\ldots,a_n)\in \mathbb{R}^n$ with $\|a\|=1$, it holds that
\begin{eqnarray*}
P\left(\left|a\cdot\epsilon\right|\leq 1\right)\geq \frac{1}{2}.
\end{eqnarray*}
This means that at least half of the $2^n$ sign sums of the form $S=\pm a_1\pm\cdots\pm a_n$ satisfy $|S|\leq 1$. In this subsection, we present a strengthened version of  Tomaszewski's conjecture, which provides information about the other half of the $2^n$ sign sums.

\begin{thm}\label{thm-2.6}
Tomaszewski's conjecture holds if and only if for any $n\in \mathbb{N}$ and   $a=(a_1,\dots,a_n)\in \mathbb{R}^n$  with $\|a\|=1$, the $2^n$ sign sums of the form $S=\pm a_1\pm\cdots\pm a_n$ can be partitioned into $2^{n-1}$ pairs $\{(s^j_1,s^j_2):j=1,\dots,2^{n-1}\}$ such  that $|s^j_1|\cdot|s^j_2|\le 1$, $j=1,\dots,2^{n-1}$.
\end{thm}

Before proving Theorem \ref{thm-2.6}, we present a lemma.

\begin{lem}\label{lem-22}
Tomaszewski's conjecture is equivalent to each of the following statements:

\noindent
(i)  For any $n\in\mathbb{N}$ and  $a=(a_1,\cdots,a_n)\in\mathbb{R}^n$,
\begin{equation}\label{eqn3-3}
    P\left(\left|a\cdot\epsilon\right|<\|a\|\right)\geq
    P\left(\left|a\cdot\epsilon\right|>\|a\|\right).
\end{equation}
(ii) For any $n\in\mathbb{N}$,
$a=(a_1,\cdots,a_n)\in\mathbb{R}^n$ with $\|a\|=1$, and $\delta>0$,
\begin{eqnarray}\label{main-inequ-a}
P\left(a\cdot\epsilon>\delta\right)
+P\left(a\cdot\epsilon>\frac{1}{\delta}\right)\leq \frac{1}{2}.
\end{eqnarray}
\end{lem}
{\bf Proof.}\ (\ref{eqn3-3}) $\Rightarrow$ (\ref{eqn3-1}):
 (\ref{eqn3-1}) follows directly from (\ref{eqn3-3}) by
\begin{eqnarray*}
P\left(|a\cdot\epsilon|\leq\|a\|\right)\geq
P\left(|a\cdot\epsilon|<\|a\|\right)\geq
P\left(|a\cdot\epsilon|>\|a\|\right)=1-P\left(|a\cdot\epsilon|\leq\|a\|\right).
\end{eqnarray*}

(\ref{eqn3-1}) $\Rightarrow$ (\ref{main-inequ-a}):
Set $2b=\delta-\frac{1}{\delta}$, $a^*=(a_1,\dots,a_n,b)$,
 $\epsilon^*=(\epsilon_1,\dots,\epsilon_n,\epsilon_{n+1})$.
Since  (\ref{eqn3-1}) holds with $(a, \epsilon)$ replaced by  $(a^*,\epsilon^*)$,
it follows that
\begin{equation}\label{eqn310}
    P\left(a^*\cdot\epsilon^*>\|a^*\|\right)\leq \frac{1}{4}.
\end{equation}
Note that
$$
a^*\cdot\epsilon^*=a\cdot\epsilon+b\,\epsilon_{n+1}, \quad
\|a^*\|=\sqrt{\|a\|^2+b^2}=\sqrt{1+b^2}. 
$$
By taking conditional probability on $\epsilon_{n+1}=\pm1$,
(\ref{eqn310}) can be written as
$$
\frac{1}{2}P\left(a\cdot\epsilon+b>\sqrt{1+b^2}\right)
+\frac{1}{2}P\left(a\cdot\epsilon-b>\sqrt{1+b^2}\right)\leq \frac{1}{4},
$$
i.e.,
$$
P\left(a\cdot\epsilon>\sqrt{1+b^2}-b\right)
+P\left(a\cdot\epsilon>\sqrt{1+b^2}+b\right)\leq \frac{1}{2}.
$$

Note that $2b=\delta-\frac{1}{\delta}$ implies that $2\sqrt{1+b^2}=\delta+\frac{1}{\delta}$ and
$$
\sqrt{1+b^2}-b=\frac{1}{\delta},\quad \sqrt{1+b^2}+b=\delta.
$$
Thus
$$
P\left(a\cdot\epsilon>\delta\right)
+P\left(a\cdot\epsilon>\frac{1}{\delta}\right)\leq \frac{1}{2},
$$
which implies that (\ref{main-inequ-a}) holds.

(\ref{main-inequ-a}) $\Rightarrow$ (\ref{eqn3-3}): By symmetry, the case that $\|a\|=0$ is trivial.

Let $\|a\|>0$.
Applying  (\ref{main-inequ-a}) to $a'=\frac{a}{\|a\|}$, we get
\begin{eqnarray*}
P\left(|a\cdot\epsilon|>\delta\|a\|\right)
+P\left(|a\cdot\epsilon|>\frac{1}{\delta}\|a\|\right)\leq 1.
\end{eqnarray*}
Letting $\delta\uparrow 1$, we obtain that
\begin{eqnarray*}
P\left(|a\cdot\epsilon|\geq\|a\|\right)
+P\left(|a\cdot\epsilon|>\|a\|\right)\leq 1.
\end{eqnarray*}
Thus
$P\left(|a\cdot\epsilon|<\|a\|\right)=1-P\left(a\cdot\epsilon\geq\|a\|\right)\geq P\left(|a\cdot\epsilon|>\|a\|\right)$, i.e., (\ref{eqn3-3}) holds. \hfill $\square$

\begin{rem}\label{rem-31}
(a) We have the following two apparently equivalent formulations of (\ref{eqn3-3}):
$$
P\left(|a\cdot\epsilon|>\|a\|\right)\leq \frac{1}{2}-\frac{1}{2}P\left(|a\cdot\epsilon|=\|a\|\right),
$$
and
$$
P(|a\cdot\epsilon|\leq \|a\|)\geq \frac{1}{2}+\frac{1}{2}P\left(|a\cdot\epsilon|=\|a\|\right).
$$

(b) An equivalent version of (\ref{main-inequ-a}) has appeared in
Dzindzalieta \cite[Proposition 1]{Dz14b} (cf. also
\cite[Proposition 13]{Dz14}):
\begin{eqnarray}\label{main-inequ-a1}
P\left(|a\cdot\epsilon|\leq\delta\right)
\geq P\left(|a\cdot\epsilon|\geq\frac{1}{\delta}\right).
\end{eqnarray}
It is shown in \cite[Proposition 1]{Dz14b} that (\ref{eqn3-1})
holds if and only if  (\ref{main-inequ-a1}) holds for all $0\le
\delta\le 1$. However, there is an error in the proof of
\cite[Proposition 1]{Dz14b}. To show that (\ref{eqn3-1}) holds
implies that (\ref{main-inequ-a1}) holds for all $0\le \delta\le 1$,
the author used the following inequality
\begin{eqnarray}\label{wrong}
P\left(|a\cdot\epsilon|\ge\delta\right)+P\left(|a\cdot\epsilon|\geq\frac{1}{\delta}\right)\le
1.
\end{eqnarray}
We point out that  (\ref{wrong}) does not
hold in general. If we let
$a=(\frac{1}{2},\frac{1}{2},\frac{1}{2},\frac{1}{2})$ and
$\delta=1$, then
$$
P\left(|a\cdot\epsilon|\ge\delta\right)+P\left(|a\cdot\epsilon|\geq\frac{1}{\delta}\right)=\frac{5}{4}>1.
$$

Although (\ref{wrong}) does not hold in general,
(\ref{main-inequ-a1}) is still an equivalent formulation of
(\ref{main-inequ-a}). In fact, suppose that (\ref{main-inequ-a})
holds for any $\delta>0$. Then, we have
\begin{eqnarray*}
P\left(|a\cdot\epsilon|\leq\delta\right)
&=&\lim_{\varepsilon\rightarrow0}P\left(|a\cdot\epsilon|\le\delta+\varepsilon\right)\\
&=&\lim_{\varepsilon\rightarrow0}\left[1-2P\left(a\cdot\epsilon>\delta+\varepsilon\right)\right]\\
&\ge&\limsup_{\varepsilon\rightarrow0}\left[1-2\left(\frac{1}{2}-P\left(a\cdot\epsilon>\frac{1}{\delta+\varepsilon}\right)\right)\right]\\
&=&\limsup_{\varepsilon\rightarrow0}P\left(|a\cdot\epsilon|>\frac{1}{\delta+\varepsilon}\right)\\
&=&P\left(|a\cdot\epsilon|\geq\frac{1}{\delta}\right),
\end{eqnarray*}
i.e., (\ref{main-inequ-a1}) holds. If (\ref{main-inequ-a1}) holds, then by the symmetry of $(\epsilon_i)$, we get
\begin{eqnarray*}
&&P\left(a\cdot\epsilon>\delta\right)
+P\left(a\cdot\epsilon>\frac{1}{\delta}\right)\\
&&=\frac{1}{2}\left[P\left(|a\cdot\epsilon|>\delta\right)
+P\left(|a\cdot\epsilon|>\frac{1}{\delta}\right)\right]\\
&&=\frac{1}{2}\left[1-P\left(|a\cdot\epsilon|\leq \delta\right)
+P\left(|a\cdot\epsilon|>\frac{1}{\delta}\right)\right]\\
&&\leq \frac{1}{2}\left[1-P\left(|a\cdot\epsilon|\leq \delta\right)
+P\left(|a\cdot\epsilon|\geq \frac{1}{\delta}\right)\right]\\
&&=\frac{1}{2}-\frac{1}{2}\left[P\left(|a\cdot\epsilon|\leq \delta\right)-P\left(|a\cdot\epsilon|\geq \frac{1}{\delta}\right)\right]\\
&&\leq \frac{1}{2},
\end{eqnarray*}
i.e., (\ref{main-inequ-a}) holds.
\end{rem}

\noindent {\bf Proof of Theorem \ref{thm-2.6}.} We need only prove the necessity. Suppose that Tomaszewski's conjecture holds. Then, Lemma \ref{lem-22} implies that (\ref{main-inequ-a}) holds.

For $a=(a_1,\dots,a_n)\in\mathbb{R}^n$ with $\|a\|=1$, we arrange the $2^n$ sign sums of the form $S=\pm a_1\pm\cdots\pm a_n$ in non-decreasing order and denote the obtained sequence by $s_{-2^{n-1}},\dots,s_{-1}, s_{1},\dots,s_{2^{n-1}}$.
It follows that $s_{-k}\leq0\leq s_{k}$ and $s_{-k}=-s_{k}$ for $k=1,\dots,2^{n-1}$. To finish the proof, it is sufficient to show that
\begin{equation}\label{eqn3-16}
 s_{k}\cdot s_{2^{n-1}+1-k}\leq1,\quad \forall k\in \{1,\dots,2^{n-1}\}.
\end{equation}

We use proof by contradiction. Assume that (\ref{eqn3-16}) does not hold.
Then, there exists some $k\in\{1,\dots,2^{n-1}\}$ such that $$s_{k}\cdot s_{2^{n-1}+1-k}>1.$$
It follows that $s_{k}>0$ and there exists a positive number
$\delta\in({1}/s_{k}, s_{2^{n-1}+1-k})$.
Thus
$$
\frac{1}{\delta}<{s_{k}},\quad\delta<s_{2^{n-1}+1-k}.
$$
Note that $\epsilon$ has a uniform distribution on $\{-1,1\}^n$ and $s_1\le\cdots\le s_{2^{n-1}}$. We get
$$
P\left(a\cdot\epsilon>\frac{1}{\delta}\right)\geq\frac{2^{n-1}+1-k}{2^n},\quad
P\left(a\cdot\epsilon>\delta\right)\geq\frac{k}{2^n}.
$$
Hence
$$
P\left(a\cdot\epsilon>\delta\right)
+P\left(a\cdot\epsilon>\frac{1}{\delta}\right)\geq\frac{2^{n-1}+1}{2^n}>\frac12,
$$
which contradicts  with (\ref{main-inequ-a}). Hence  (\ref{eqn3-16}) holds and therefore the proof is complete. \hfill\fbox

\section{On Hitczenko-Kwapie\'{n}'s conjecture}\setcounter{equation}{0}

In this section, we first consider Hitczenko-Kwapie\'{n}'s conjecture and prove the following result.

\begin{thm}\label{thm-21}
Suppose $n\leq7$.
Then for any $a\in \mathbb{R}^n$, it holds that
$$
     P\left(|a\cdot\epsilon|\geq\|a\|\right)\geq\frac{7}{32}.
$$
\end{thm}

\begin{rem}\label{rem-21}
Since $P\left(|a\cdot\epsilon|\geq\|a\|\right)=\frac{7}{32}$ for $a=(1,1,1,1,1,1,0)$,
the lower bound $\frac{7}{32}$ is optimal. Let $G_n$  be defined by (\ref{Gn}). Then, Theorem \ref{thm-21} implies that
$G_7=\frac{7}{32}$.
\end{rem}

We will also consider the constant $G_n'$ defined by (\ref{Gn'}) and prove the following result.

\begin{thm}\label{thm-22}
$
G_1'=0,\ G_2'=\frac{1}{2},\ G_3'=G_5'=\frac{1}{4},\ G_4'=\frac{1}{8},\ G_6'=\frac{3}{16},\ G_7'=\frac{7}{32}.
$
\end{thm}

\subsection{Some lemmas}

By the symmetry of $(\epsilon_i)$, we can assume without loss of generality that $a_1\geq a_2\geq \cdots\geq a_n\geq 0$.
Define
$$
R^{(n)}_{d}=\{(a_1,\dots,a_n)\in \mathbb{R}^n:a_1\geq a_2\geq \cdots\geq a_n\geq 0\}.
$$
For any sign vector $s=(s_1,\dots,s_n)\in\{-1,1\}^n$,
set $a\cdot{s}=a_1{s}_1+\cdots+a_n{s}_n.$
We will estimate the cardinality of the set
$$
V_{sd}(a)=\{s\in\{-1,1\}^n: a\cdot{s}\geq\|a\|\}.
$$

For two different vectors $s$ and $s'$ in $\{-1,1\}^n$,
we say that $s'$ is better than $s$ if $a\cdot s'\geq a\cdot s$
for all $a\in R^{(n)}_{d}$ and denote it by $s\prec s'$.
If $s\prec s'$ or $s=s'$, we denote $s\preceq s'$. It is easy to see that
if $s\in V_{sd}(a)$ and $s'$ is better than $s$, then $s'$ is also in $V_{sd}(a)$.
Our idea is to find an $s$ in $V_{sd}(a)$ such that
there are as many as possible sign vectors which are better than $s$.

For ${J}\subseteq I_n=\{1,\dots,n\}$,
let $({J})_n$ denote the sign vector $s=(s_1,\dots,s_n)\in\{-1,1\}^n$
such that $s_j=-1$ for $j\in {J}$ and $s_j=1$ otherwise.
Obviously ${J}\mapsto ({J})_n$ sets up a one-to-one correspondence between $2^{I_n}$ and $\{-1,1\}^n$.
For example,
$(6)_7$ and $(3,4)_7$ denote $(1,1,1,1,1,-1,1)$ and $(1,1,-1,-1,1,1,1)$, respectively.
In the case that ${J}=\varnothing$, we just write $( )_7$ for $(1,1,1,1,1,1,1)$.
The following lemma can be easily checked.
\begin{lem}\label{pro-21} Suppose that $J,J',J''\subseteq  I_n$.

\noindent (i) If $i<j$, then $(i)_n\prec (j)_n$.

\noindent (ii) If $J\supset J'$, then $({J})_n\prec ({J}')_n$.

\noindent (iii) If $({J})_n\prec ({J}')_n$ and $({J}')_n\prec ({J}'')_n$,
then $({J})_n\prec (J'')_n$.

\noindent (iv) If $({J}')_n\prec ({J}'')_n$ and ${J}\cap J'={J}\cap J''=\varnothing$,
then $({J}\cup J')_n\prec ({J}\cup J'')_n$.

\noindent (v) Suppose that $ i_1\leq j_1,\dots,i_m\leq j_m$ for some
$i_1<i_2<\cdots<i_m$ and $j_1<j_2<\dots<j_m$.
Then $$(i_1,\dots,i_m)_n\preceq (j_1,\dots,j_m)_n.$$
\end{lem}

For ${J}\subseteq I_n$, let $J^{(2)}$ be defined by (\ref{JJ}).

\begin{lem}\label{lem220}
Let $a\in R^{(n)}_{d}$. Suppose that
${J},K\subseteq I_n$ with $a\cdot({J})_n\geq 0$ and $a\cdot({K})_n\geq 0$.

\noindent (i) Let  ${J}_c=I_n-{J}$. Then $({J})_n\in V_{sd}(a)$ if and only if
$$
    \sum_{(i,j)\in J^{(2)}\cup J_c^{(2)}} {a_i}{a_j}-
    \sum_{(i, j)\in J\times J_c} a_i a_j\geq 0.
$$
(ii) Let ${({J}\cup{K})}_c=I_n-({J}\cup{K})$.
Suppose that ${J}\cap{K}=\varnothing$ and
\begin{equation}\label{eqn222}
    \sum_{(i,j)\in J^{(2)}\cup K^{(2)}\cup {({J}\cup{K})}_c^{(2)}} {a_i}{a_j}
    -\sum_{(i, j)\in J\times K} a_i a_j\geq 0.
\end{equation}
Then, either $({J})_n$ or $({K})_n$ is in  $V_{sd}(a)$.
\end{lem}

\noindent {\bf Proof.} (i) Note that $a\cdot({J})_n=\sum_{i\in J_c}a_i-\sum_{i\in J}a_i$. By the assumption that $a\cdot({J})_n\geq 0$, we get
\begin{eqnarray*}
({J})_n\in V_{sd}(a)
&\Leftrightarrow&\left[\, a\cdot({J})_n\, \right]^2\geq \|a\|^2\nonumber\\
&\Leftrightarrow&\left[\sum_{i\in J_c}a_i-\sum_{i\in J}a_i\right]^2\geq\sum_{i\in I_n}a_i^2\\
&\Leftrightarrow&\sum_{(i,j)\in J^{(2)}\cup J_c^{(2)}} {a_i}{a_j}-
    \sum_{(i, j)\in J\times J_c} a_i a_j\geq 0.
\end{eqnarray*}

(ii) Note that ${J}_c={({J}\cup{K})}_c\cup K$ and ${K}_c={({J}\cup{K})}_c\cup J$. Then
\begin{eqnarray*}
&&\sum_{(i,j)\in J^{(2)}\cup J_c^{(2)}} {a_i}{a_j}
-\sum_{(i, j)\in K\times K_c} a_i a_j \\
&&=
\left[\sum_{(i,j)\in J^{(2)}\cup K^{(2)}\cup {({J}\cup{K})}_c^{(2)}} {a_i}{a_j}
+\sum_{(i,j)\in K\times {({J}\cup{K})}_c} {a_i}{a_j}\right]
-\left[\sum_{(i,j)\in K\times {({J}\cup{K})}_c} {a_i}{a_j}
+\sum_{(i, j)\in K\times J} a_i a_j \right]\\
&&= \sum_{(i,j)\in J^{(2)}\cup K^{(2)}\cup {({J}\cup{K})}_c^{(2)}} {a_i}{a_j}
    -\sum_{(i, j)\in J\times K} a_i a_j.
\end{eqnarray*}
Thus, by symmetry,
\begin{eqnarray*}
&&
\left[\sum_{(i,j)\in J^{(2)}\cup J_c^{(2)}} {a_i}{a_j}-\sum_{(i, j)\in J\times J_c} a_i a_j\right]
+\left[\sum_{(i,j)\in K^{(2)}\cup K_c^{(2)}} {a_i}{a_j}-\sum_{(i, j)\in K\times K_c} a_i a_j\right] \\
&&= 2\left[\sum_{(i,j)\in J^{(2)}\cup K^{(2)}\cup {({J}\cup{K})}_c^{(2)}} {a_i}{a_j}
    -\sum_{(i, j)\in J\times K} a_i a_j \right].
\end{eqnarray*}
By (\ref{eqn222}), we find that at least one term of the left hand side is non-negative.
Therefore, either $({J})_n$ or $({K})_n$ is in  $V_{sd}(a)$ by (i). \hfill $\square$

\begin{lem}\label{lem222}
Let $a\in R^{(7)}_{d}$. Then at least one of  $(2)_7, (3,4)_7$ and  $(5,6,7)_7$ is in $V_{sd}(a)$.
\end{lem}

\noindent {\bf Proof.} {\it Case 1}: $a_2< a_3+a_5+a_6+a_7$.

Applying Lemma \ref{lem220} with $J=\{2\}$ and $K=\{5,6,7\}$,
${({J}\cup{K})}_c=\{1,3,4\}$, we get
\begin{eqnarray*}
&&(a_1a_3+a_1a_4+a_3a_4)+(a_5a_6+a_5a_7+a_6a_7)-a_2(a_5+a_6+a_7)\\
&&=(a_1a_3+a_1a_4-a_2a_5-a_2a_6)+(a_3a_4+a_5a_6+a_5a_7+a_6a_7)-a_2a_7\\
&&\geq (a_3a_4+a_5a_6+a_5a_7+a_6a_7)-(a_3+a_5+a_6+a_7)a_7\\
&&\geq 0.
\end{eqnarray*}
Then, we obtain by Lemma \ref{lem220} (ii) that either $(2)_7$ or $(5,6,7)_7$ is in  $V_{sd}(a)$.

{\it Case 2}: $a_2\geq a_3+a_5+a_6+a_7$.

Applying Lemma \ref{lem220} with $J=\{3,4\}$ and $K=\{5,6,7\}$, we get
\begin{eqnarray*}
&&a_1a_2+a_3a_4+(a_5a_6+a_5a_7+a_6a_7)-(a_3+a_4)(a_5+a_6+a_7)\\
&&\geq a_2^2-2a_3(a_5+a_6+a_7)\\
&&\geq (a_3+a_5+a_6+a_7)^2-2a_3(a_5+a_6+a_7)\\
&&\geq 0.
\end{eqnarray*}
Then, we obtain by Lemma \ref{lem220} (ii) that either $(3,4)_7$ or $(5,6,7)_7$ is in  $V_{sd}(a)$.

Therefore, for both case 1 and case 2,
at least one of  $(2)_7$, $(3,4)_7$, $(5,6,7)_7$ is in  $V_{sd}(a)$. \hfill $\square$

\subsection{Proof of Theorem \ref{thm-21}}
Since $\epsilon=(\epsilon_1,\dots,\epsilon_7)$
has a uniform distribution on $\{-1,1\}^7$, which consists of $2^7=128$ sign vectors. We have
$$
P\left(a\cdot\epsilon\geq\|a\|\right)=\frac{|V_{sd}(a)|}{128}.
$$
Then it is sufficient to show that $V_{sd}(a)$ contains at least 14
elements for any $a\in R^{(7)}_{d}$.

\textbf{Case 1}:  $a_1\geq a_4+a_5$.

Applying Lemma \ref{lem220} with $J=\{3,6\}$ and $K=\{4,5,7\}$,
${({J}\cup{K})}_c=\{1,2\}$,  we get
\begin{eqnarray*}
&&a_1a_2+a_3a_6+(a_4a_5+a_4a_7+a_5a_7)-(a_3+a_6)(a_4+a_5+a_7)\\
&&=a_1a_2-a_3(a_4+a_5)+(a_3-a_5)(a_6-a_7)+a_4(a_5-a_6)+(a_4-a_6)a_7\\
&&\geq 0.
\end{eqnarray*}
Then, we obtain by Lemma \ref{lem220} (ii) that either
$(3,6)_7$ or $(4,5,7)_7$ is in $V_{sd}(a)$.

Case 1a: Suppose that $(4,5,7)_7\in V_{sd}(a)$.
By Lemma \ref{pro-21}, the following 13 sign vectors are all better than $(4,5,7)_7$:
$$
(4,6,7)_7, \ (5,6,7)_7,\ (4,5)_7, \ (4,6)_7, \ (4,7)_7, \ (5,6)_7, \  (5,7)_7, \ (6,7)_7, \
(4)_7, \ (5)_7, \ (6)_7, \ (7)_7, \ ()_7.
$$
 Then all of these vectors are in $V_{sd}(a)$.
 Hence there are at least 14 elements in $V_{sd}(a)$.

Case 1b: Suppose that $(3,6)_7\in V_{sd}(a)$.
Then the following 12 better sign vectors are also in $V_{sd}(a)$:
$$
(3,7)_7,\ (4,6)_7, \ (4,7)_7, \ (5,6)_7, \  (5,7)_7, \ (6,7)_7, \
(3)_7, \ (4)_7, \ (5)_7, \ (6)_7, \ (7)_7, \ ()_7.
$$
Moreover, by Lemma \ref{lem222},
 at least one of  $(2)_7$, $(3,4)_7$, $(5,6,7)_7$ is in  $V_{sd}(a)$,
 and the three sign vectors are not included in the above sequence.
 Thus there are also at least 14 elements in $V_{sd}(a)$.

\textbf{Case 2}:  $a_1< a_4+a_5$.

We will show that in this case
 at least one of  $(1,7)_7$, $(3,6)_7$, $(4,5)_7$ is in  $V_{sd}(a)$.
 Since $a\cdot(3,6)_7$, $a\cdot(4,5)_7$ and $a\cdot(1,7)_7$ are all nonnegative,
 it is sufficient to show that
 $$
 \left[\,a\cdot(3,6)_7\right]^2+\left[\,a\cdot(4,5)_7\right]^2+\left[\,a\cdot(1,7)_7\right]^2
 -3 \|a\|^2\geq0.
 $$
 Expanding the left hand side of the above inequality, we get
\begin{eqnarray*}
&&\left[\,a\cdot(3,6)_7\right]^2+\left[\,a\cdot(4,5)_7\right]^2+\left[\,a\cdot(1,7)_7\right]^2
 -3 \|a\|^2\\
&&=3(a_1a_7+a_3a_6+a_4a_5)+a_2(a_1+a_3+a_4+a_5+a_6+a_7)\nonumber\\
&&\quad-(a_1+a_7)(a_3+a_6)-(a_1+a_7)(a_4+a_5)-(a_3+a_6)(a_4+a_5)\nonumber\\
&&=-a_1(a_4+a_5+a_6)+(a_2-a_3)(a_1+a_4+a_5)+(2a_3-a_4-a_5)a_6\nonumber\\
&&\quad+(3a_1+a_2-a_3-a_4-a_5-a_6)a_7+a_2a_3+a_2a_6+a_3a_6+3a_4a_5\nonumber\\
&&\geq a_2a_3+a_2a_6+a_3a_6+3a_4a_5-a_1(a_4+a_5+a_6)\nonumber\\
&&> a_2a_3+(a_2+a_3)a_6+3a_4a_5-(a_4+a_5)(a_4+a_5+a_6)\nonumber\\
&&\geq 0.
\end{eqnarray*}
Then it follows that at least one of $(1,7)_7$, $(3,6)_7$, $(4,5)_7$ is in $V_{sd}(a)$.

Case 2a: Suppose that $(1,7)_7\in V_{sd}(a)$.
Then the following 13 better sign vectors are also in $V_{sd}(a)$:
$$
 (2,7)_7, \ (3,7)_7, (4,7)_7, \ (5,7)_7,  \ (6,7)_7, \
 (1)_7, \ (2)_7, \ (3)_7, \ (4)_7, \ (5)_7, \ (6)_7, \ (7)_7, \ ()_7.
$$
Thus there are at least 14 elements in $V_{sd}(a)$.

Case 2b: Suppose that $(3,6)_7\in V_{sd}(a)$.
Then as in Case 1b there are at least 14 elements in $V_{sd}(a)$.

Case 2c: Suppose that $(4,5)_7\in V_{sd}(a)$.
Then the following 10 better sign vectors are also in $V_{sd}(a)$:
 $$
 (4,6)_7, \ (4,7)_7, \ (5,6)_7, \  (5,7)_7, \ (6,7)_7,
 \ (4)_7, \ (5)_7, \ (6)_7, \ (7)_7, \ ()_7.
 $$
 Moreover, since $a_1< a_4+a_5$ and $(4,5)_7\in V_{sd}(a)$,
 we have $$a\cdot(1)_7>a\cdot(4,5)_7\geq\|a\|.$$
 Then $(1)_7$ is in $V_{sd}(a)$. Hence $(2)_7$ and $(3)_7$ are also in $V_{sd}(a)$.
 Therefore, in this case there are at least 14 elements in $V_{sd}(a)$.

\subsection{Proof of Theorem \ref{thm-22}}
By the symmetry of $(\epsilon_i)$, we can assume without loss of generality that  $a_1\geq a_2\geq \cdots\geq a_n>0$.
Define
$$
V_{sd+}(a)=\{s\in\{-1,1\}^n: a\cdot{s}>\|a\|\}.
$$
Obviously, if $s\in V_{sd+}(a)$ and $s'$ is better than $s$, then  $s'$ is also in $V_{sd+}(a)$. Note that
$$
G'_n=2\,\inf\left\{ P\left(a\cdot\epsilon>\|a\|\right):{a\in \mathbb{R}^n}, a_i\neq0,\forall i=1,\ldots,n\right\}.
$$

(1)  By the definition of $G_n'$, we have $G_1'=0$.

(2)  By $a_1-a_2<\sqrt{a_1^2+a_2^2}<a_1+a_2$,
we have $P\left(a\cdot\epsilon>\|a\|\right)=\frac{1}{4}$ for $n=2$.
Then $G_2'=\frac{1}{2}$.

(3)  By $a_1+a_2+a_3>\sqrt{a_1^2+a_2^2+a_3^3}$, we know that $( )_3\in V_{sd+}(a)$.
Then $P\left(a\cdot\epsilon>\|a\|\right)\geq\frac{1}{8}$ for $n=3$.
On the other hand,
$P\left(a\cdot\epsilon>\|a\|\right)=\frac{1}{8}$ for $a=(1,1,1)$.
Thus $G_3'=\frac{1}{4}$.

(4)  As above, we have $( )_4\in V_{sd+}(a)$.
Then $P\left(a\cdot\epsilon>\|a\|\right)\geq\frac{1}{16}$ for $n=4$.
On the other hand,
$P\left(a\cdot\epsilon>\|a\|\right)=\frac{1}{16}$ for $a=(1,1,1,1)$.
Thus $G_4'=\frac{1}{8}$.

Before presenting the proof for $n=5,6,7$, we need the following lemma whose proof is similar to that of  Lemma \ref{lem220}.
\begin{lem}\label{lem230}
Let $a=(a_1,\dots,a_n)$ be a vector with $a_1\geq \cdots\geq a_n>0$.
Suppose that ${J},K\subseteq I_n$ satisfying ${J}\cap{K}=\varnothing$, $a\cdot({J})_n>0$, $a\cdot({K})_n>0$ and
$$
    \sum_{(i,j)\in J^{(2)}\cup K^{(2)}\cup {({J}\cup{K})}_c^{(2)}} {a_i}{a_j}
    -\sum_{(i, j)\in J\times K} a_i a_j>0.
$$
Then either $({J})_n$ or $({K})_n$ is in $V_{sd+}(a)$.
\end{lem}

(5)  Let $n=5$. We will show that $V_{sd+}(a)$ contains at least 4 elements.
Applying Lemma \ref{lem230} with $J=\{3\}$ and $K=\{4,5\}$,
${({J}\cup{K})}_c=\{1,2\}$, we get
$$
    a_1a_2+a_4a_5-a_3a_4-a_3a_5=a_1a_2-a_3^2+(a_3-a_4)(a_3-a_5)\geq0.
$$

\textbf{Case (i):}  Suppose that $a_1a_2-a_3^2+(a_3-a_4)(a_3-a_5)>0$.
Then, by Lemma \ref{lem230} (ii), we find that either $(3)_5$ or $(4,5)_5$ is in $V_{sd+}(a)$.
Note that $(4)_5$, $(5)_5$ and $()_5$ are better than both $(3)_5$ and $(4,5)_5$.
In either case $V_{sd+}(a)$ contains at least 4 elements. Thus
\begin{equation}\label{prf22p}
    P\left(a\cdot\epsilon>\|a\|\right)\geq\frac{4}{32}=\frac{1}{8}.
\end{equation}

\textbf{Case (ii):} Suppose that $a_1a_2-a_3^2+(a_3-a_4)(a_3-a_5)=0$.
Then we have $a_1=a_2=a_3=a_4\geq a_5>0$.
Thus $a\cdot (1)_5=2a_1+a_5>\sqrt{4a_1^2+a_5^2}=\|a\|$,
which implies that $(1)_5\in V_{sd+}(a)$,
and hence $(2)_5, \ (3)_5, \ (4)_5, \ (5)_5, \ ()_5$ are all in $V_{sd+}(a)$.
Therefore $P\left(a\cdot\epsilon>\|a\|\right)=\frac{6}{32}$ and (\ref{prf22p}) still holds.

On the other hand, we have
$P\left(a\cdot\epsilon>\|a\|\right)=\frac{1}{8}$ for $a=(2,2,1,1,1)$.
Hence $G_5'=\frac{1}{4}$.

(6)  Let $n=6$. We will show that $V_{sd+}(a)$ contains at least 6 elements.
Applying Lemma \ref{lem230} with $J=\{2\}$ and $K=\{4,6\}$,
${({J}\cup{K})}_c=\{1,3,5\}$, we get
\begin{eqnarray*}
 a_4a_6+(a_1a_3+a_1a_5+a_3a_5)-(a_2a_4+a_2a_6)>
 (a_1a_3-a_2a_4)+(a_1a_5-a_2a_6)\geq0.
\end{eqnarray*}
Then, by Lemma \ref{lem230} (ii), we find that either $(2)_6$ or $(4,6)_6$ is in $V_{sd+}(a)$.

\textbf{Case (i):} Suppose that $(2)_6$ is in $V_{sd+}(a)$.
Then the following better sign vectors are also in $V_{sd}(a)$:
$$
 (3)_6, \ (4)_6, \ (5)_6, \ (6)_6, \ ()_6.
$$
Thus there are at least 6 elements in $V_{sd}(a)$.

\textbf{Case (ii):} Suppose that $(4,6)_6$ is in $V_{sd+}(a)$.
Then the following sign vectors are also in $V_{sd}(a)$:
$$
  (5,6)_6,  \ (4)_6, \ (5)_6, \ (6)_6, \ ()_6.
$$
Thus in either case $V_{sd+}(a)$ contains at least 6 elements. Hence
$$P\left(a\cdot\epsilon>\|a\|\right)\geq\frac{6}{64}=\frac{3}{32}.$$
On the other hand,
$P\left(a\cdot\epsilon>\|a\|\right)=\frac{3}{32}$ for $a=(2,1,1,1,1,1)$.
Therefore $G_6'=\frac{3}{16}$.

(7)  Let $n=7$. The proof for this case is similar to that of Theorem \ref{thm-21}.
By Lemma \ref{lem230} and the proof of Lemma \ref{lem222}, we can prove the following lemma,
since the equality in that proof does hold in the case that $a_1\geq\cdots\geq a_7>0$.
\begin{lem}\label{lem237}
Suppose that $a=(a_1,\dots,a_7)$ with $a_1\geq\cdots\geq a_7>0$.
Then at least one of  $(2)_7, (3,4)_7$ and  $(5,6,7)_7$ is in $V_{sd+}(a)$.
\end{lem}

We now show that $V_{sd+}(a)$ contains at least 14 elements.

\textbf{Case 1}:  $a_1>a_4+a_5$.
Applying Lemma \ref{lem230} with $J=\{3,6\}$ and $K=\{4,5,7\}$,
${({J}\cup{K})}_c=\{1,2\}$, following the proof of Theorem \ref{thm-21} (Case 1), we get
\begin{eqnarray*}
&&a_1a_2+a_3a_6+(a_4a_5+a_4a_7+a_5a_7)-(a_3+a_6)(a_4+a_5+a_7)\\
&&= a_1a_2-a_3(a_4+a_5)+(a_3-a_5)(a_6-a_7)+a_4(a_5-a_6)+(a_4-a_6)a_7\\
&&> 0.
\end{eqnarray*}
Then, we obtain by  Lemma \ref{lem230} (ii) that either
$(3,6)_7$ or $(4,5,7)_7$ is in $V_{sd+}(a)$.

Case 1a: Suppose that $(4,5,7)_7\in V_{sd+}(a)$.
By Lemma \ref{pro-21}, the following 13 sign vectors are all better than $(4,5,7)_7$:
$$
(4,6,7)_7, \ (5,6,7)_7,\ (4,5)_7, \ (4,6)_7, \ (4,7)_7, \ (5,6)_7, \  (5,7)_7, \ (6,7)_7, \
(4)_7, \ (5)_7, \ (6)_7, \ (7)_7, \ ()_7.
$$
 Then all of these vectors are  in $V_{sd+}(a)$.
 Hence there are at least 14 elements in $V_{sd+}(a)$.

Case 1b: Suppose that $(3,6)_7\in V_{sd+}(a)$.
Then the following 12 better sign vectors are also in $V_{sd+}(a)$:
$$
(3,7)_7,\ (4,6)_7, \ (4,7)_7, \ (5,6)_7, \  (5,7)_7, \ (6,7)_7, \
(3)_7, \ (4)_7, \ (5)_7, \ (6)_7, \ (7)_7, \ ()_7.
$$
 Moreover, by Lemma \ref{lem237},
 at least one of  $(2)_7$, $(3,4)_7$, $(5,6,7)_7$ is in  $V_{sd+}(a)$,
 and the 3 sign vectors are not included in the above sequence.
 Thus there are also at least 14 elements in $V_{sd+}(a)$.

\textbf{Case 2}:  $a_1\leq a_4+a_5$.
We will show that in this case
 at least one of $(1,7)_7$, $(3,6)_7$, $(4,5)_7$ is in  $V_{sd+}(a)$.
  Since $a\cdot(3,6)_7$, $a\cdot(4,5)_7$ and $a\cdot(1,7)_7$ are all nonnegative,
 it is sufficient to show that
\begin{equation}\label{eqc7ca2}
     \left[\,a\cdot(3,6)_7\right]^2+\left[\,a\cdot(4,5)_7\right]^2
     +\left[\,a\cdot(1,7)_7\right]^2-3 \|a\|^2>0.
\end{equation}
Expanding the left hand side of the above inequality, we get
\begin{eqnarray*}
&&\left[\,a\cdot(3,6)_7\right]^2+\left[\,a\cdot(4,5)_7\right]^2
     +\left[\,a\cdot(1,7)_7\right]^2-3 \|a\|^2\\
&&=3(a_1a_7+a_3a_6+a_4a_5)+a_2(a_1+a_3+a_4+a_5+a_6+a_7)\nonumber\\
&&\quad-(a_1+a_7)(a_3+a_6)-(a_1+a_7)(a_4+a_5)-(a_3+a_6)(a_4+a_5)\nonumber\\
&&=-a_1(a_4+a_5+a_6)+(a_2-a_3)(a_1+a_4+a_5)+(2a_3-a_4-a_5)a_6\nonumber\\
&&\quad+(3a_1+a_2-a_3-a_4-a_5-a_6)a_7+a_2a_3+a_2a_6+a_3a_6+3a_4a_5\nonumber\\
&&\geq a_2a_3+a_2a_6+a_3a_6+3a_4a_5-a_1(a_4+a_5+a_6)\nonumber\\
&&\geq a_2a_3+(a_2+a_3)a_6+3a_4a_5-(a_4+a_5)(a_4+a_5+a_6)\nonumber\\
&&\geq  0.
\end{eqnarray*}

Note that the equalities in the above three ``$\geq$'' signs cannot hold simultaneously.
In fact, the equality in the first ``$\geq$'' sign implies that $a_1=a_2=a_3=a_4=a_5=a_6$,
while the equality in the second ``$\geq$'' sign implies that $a_1=a_4+a_5$.
Then the equalities hold simultaneously if and only if  $a_1=0$
(and hence $a_2=\cdots=a_7=0$),
which is in conflict with the assumption of Theorem \ref{thm-22}.
Thus the equalities cannot hold simultaneously, i.e., (\ref{eqc7ca2}) holds.
Therefore, it follows that at least one of $(1,7)_7$, $(3,6)_7$, $(4,5)_7$ is in $V_{sd+}(a)$.

Case 2a: Suppose that $(1,7)_7\in V_{sd+}(a)$.
Then the following 13 sign vectors are also in $V_{sd+}(a)$:
$$
 (2,7)_7, \ (3,7)_7, (4,7)_7, \ (5,7)_7,  \ (6,7)_7, \
 (1)_7, \ (2)_7, \ (3)_7, \ (4)_7, \ (5)_7, \ (6)_7, \ (7)_7, \ ()_7.
$$
Thus there are at least 14 elements in $V_{sd+}(a)$.

Case 2b: Suppose that $(3,6)_7\in V_{sd+}(a)$.
Then as in Case 1b there are at least 14 elements in $V_{sd+}(a)$.

Case 2c: Suppose that $(4,5)_7\in V_{sd+}(a)$.
Then the following 10 sign vectors are also in $V_{sd+}(a)$:
 $$
 (4,6)_7, \ (4,7)_7, \ (5,6)_7, \  (5,7)_7, \ (6,7)_7,
 \ (4)_7, \ (5)_7, \ (6)_7, \ (7)_7, \ ()_7.
 $$
 Moreover, since $a_1\leq a_4+a_5$ and $(4,5)_7\in V_{sd+}(a)$,
 we have $$a\cdot(1)_7\geq a\cdot(4,5)_7>\|a\|.$$
 Then $(1)_7$ is in $V_{sd+}(a)$. Thus $(2)_7$ and $(3)_7$ are also in $V_{sd+}(a)$.
 Therefore, in this case there are also at least 14 elements in $V_{sd+}(a)$.

Now we have shown that  in either case  $V_{sd+}(a)$ contains at least 14 elements. Then
$$P\left(a\cdot\epsilon>\|a\|\right)\geq\frac{14}{2^7}=\frac{7}{64}.$$
On the other hand,
$P\left(a\cdot\epsilon>\|a\|\right)=\frac{7}{64}$ for $a=(2,2,2,1,1,1,1)$.
Therefore $G_7'=\frac{7}{32}$. \hfill $\square$

\bigskip

{ \noindent {\bf\large Acknowledgments}\ \  This work was supported by the National Natural Science Foundation of China (No. 11771309 and No. 11871184), the China Scholarship Council (No. 201809945013)  and the Natural Sciences and Engineering Research Council of Canada.


\begin{thebibliography}{1234}

\bibitem{BNR02} A. Ben-Tal, A. Nemirovski and C. Roos. Robust solutions of uncetain quadratic and conic-quadratic problems. {\it SIAM J. Optim.}, 13: 535-560, 2002.

\bibitem{BH17} R. B. Boppana and R. Holzman. Tomaszewsk's problem on randomly signed sums: breaking the $3/8$ barrier. arXiv: 1704.00350v3, 2017.

\bibitem{Bu67} D. L. Burkholder. Independent sequences with the Stein property. {\it Ann. Math. Statis.}, 39: 1282-1288, 1967.

\bibitem{DDS16} A. De, I. Diakonikolas and R. A. Servedio. A robust Khintchine inequality, and algorithms for computing optimal constants in Fourier analysis and high-dimensional geoemtry. {\it SIAM J. Disc. Math.}, 30: 1058-1094, 2016.

\bibitem{Dz14b} D. Dzindzalieta.    A note on random signs. {\it Lithuanian Math. J.}, 54: 403-408, 2014.

\bibitem{Dz14} D. Dzindzalieta. Tight Bernoullie tail probability bounds. PhD thesis, Vilnius University, 2014.

\bibitem{Gu86} R. K. Guy. Any answers anent these analytical enigmas?. {\it Amer. Math. Mon.}, 93: 279-281, 1986.

\bibitem{HZ17} H. Hendriks and C. A. van Zuijlen. Linear combinatorics of Rademacher random variables. arXiv: 1703.07251v1, 2017.

\bibitem{He12} F. von Heymann. Ideas for an old analytic enigma about the sphere that fail in intriguing ways. http://www.uni-koeln.de/opt/wp-content/uploads/2017/02/Cubes\_phere.pdf, 2012.

\bibitem{HK94} P. Hitczenko and S. Kwapie\'{n}. On the Rademacher series. In {\it Probability in Banach spaces, 9 (Sandjberg, 1994)}, Volume 35 of {\it Progr. Probab.}, 31-36. Birkh\"{a}user Boston, Boston, MA, 1994.

\bibitem{HK92} R. Holzman and D. J. Kleiman. On the product of sign vectors and unit vectors. {\it Combinatorica}, 12: 303-316, 1992.



\bibitem{Ol96} K. Oleszkiewicz. On the Stein property of Rademacher sequences. {\it Probab. Math. Stat.}, 16: 127-130, 1996.

\bibitem{Sh12} I. Shnurnikov. On a sum of centered random variables with nonreducing variances. arXiv:1202.2990v2, 2012.

\bibitem{To18} T. Toufar. Tomaszewski's conjecture. Master thesis, Charles University, 2018.

\bibitem{Zu11} M. C. A. van Zruilen. On a conjecture concerning the sum of independent Rademacher random variables. arXiv:1112.4988v1, 2011.


\end{thebibliography}
\end{document}